\documentclass[10pt]{amsart}
\usepackage{amssymb}
\usepackage{amscd}
\usepackage{color}

\theoremstyle{definition}
\newtheorem{thm}{Theorem}[section]

\newtheorem{prp}[thm]{Proposition}
\newtheorem{dfn}[thm]{Definition}
\newtheorem{cor}[thm]{Corollary}

\newtheorem{rmk}[thm]{Remark}

\newcommand{\beq}{\begin{equation}}
\newcommand{\eeq}{\end{equation}}
\newcommand{\beqr}{\begin{eqnarray*}}
\newcommand{\eeqr}{\end{eqnarray*}}
\newcommand{\bal}{\begin{align*}}
\newcommand{\eal}{\end{align*}}
\newcommand{\bei}{\begin{itemize}}
\newcommand{\eei}{\end{itemize}}

\newcommand{\af}{\alpha}

\newcommand{\N}{{\mathbf{N}}}

\newcommand{\tsr}{{\mathrm{tsr}}}

\newcommand{\ca}{C*-algebra}

\newcommand{\Aut}{{\rm Aut}}

\title{Stable rank of inclusion of C*-algebras of depth 2}

\author{Hiroyuki Osaka$^{*}$}

\thanks{$*$ Research partially supported by
Open Research Center Project for Private Universities : matching fund
from MEXT, 2004 - 2008 and the Grant-in Aid for Scientific Research,
Ritsumeikan University, 2005(H. Osaka).}
\address{ Department of Mathematical Sciences\\
 Ritsumeikan University\\ Kusatsu, Shiga, 520 - 2152   Japan}

\email[]{osaka@se.ritsumei.ac.jp}

\author{Tamotsu Teruya$^{**}$}

\address{ Department of Mathematical Sciences\\
 Ritsumeikan University\\ Kusatsu, Shiga, 520 - 2152   Japan}

\email[]{teruya@se.ritsumei.ac.jp}

\date{}

\begin{document}

\maketitle

\begin{abstract}
Let $1 \in A \subset B$ be an inclusion of unital C*-algebras 
of index-finite type and depth $2$.
Suppose that $A$ is infinite dimensional simple with $\tsr(A) = 1$ 
and SP-property.
Then $\tsr(B) \leq 2$. 
As a corollary
when $A$ is a simple \ca~ 
with $\tsr(A) = 1$ and SP-property 
and $\af$ an action of a finite group $G$ on $\Aut(A)$, 
$
\tsr(A \rtimes_\af G) \leq 2.
$
\end{abstract}

\section{Introduction}
The notion of topological stable rank for a C*-algebra $A$,
 denoted by $\tsr(A)$, was introduced by Rieffel, which
generalizes the concept of dimension of a topological space
(\cite{Rf1}).
He presented the basic properties and stability theorem related to
K-Theory for C*-algebras. In \cite{Rf1} he proved that
$\tsr(A \rtimes_\alpha {\mathbb Z}) \leq \tsr(A) + 1$, and asked if
an irrational rotation algebra $A_\theta$ has topological stable rank two.
I. Putnum (\cite{pu}) gave a complete answer to this question, that is,
$\tsr(A_\theta) = 1$. Moreover, using the notion of approximate
divisibility and U. Haggerup's striking result (\cite{ha}, \cite{ht}),
Blackadar, Kumjian, and R\o rdam (\cite{bkr}) proved that
every nonrational
 noncommutative torus has topological stable rank one.
Naturally, we pose a question of how to compute topological
stable rank of $A \rtimes_\alpha G$ for any discrete group $G$.

On the contrary, one of long standing problems was whether the fixed point
algebra of a UHF C*-algebra by an action of a finite group $G$ is an
AF C*-algebra. In 1988, Blackadar (\cite{bl3}) constructed a
symmetry
on the CAR algebra whose fixed point algebra is not an AF C*-algebra.
Note that Kumjian (\cite{km}) constructed a symmetry
on a simple AF C*-algebra
whose fixed point algebra is not an AF C*-algebra.
Blackadar proposed
the question in \cite{bl3} whether $\tsr(A \rtimes_\alpha G) = 1$
for any unital AF C*-algebra $A$,
a finite group $G$, and an action $\alpha$ of $G$ on $A$.

In \cite{OT} the authors presented a partial answer to an 
extended question of Blackadr's using C*-index Theory by Watatani (\cite{wata}), 
that is, 
Let $1 \in A \subset B$ be an inclusion of unital
C*-algebras and $E\colon B \rightarrow A$ be
a faithful conditional expectation
of index-finite type.
Suppose that the inclusion $1 \in A \subset B$
has depth $2$ and
$A$ is tsr boundedly divisible with $\tsr(A) = 1$.
Then $\tsr(B) \leq 2$. 
Here a C*-algebra $A$ is {\it tsr boundedly divisible}
(\cite[Definition 4.1]{rf3})
if there is a constant $K$ ($> 0$) such that for every positive integer $m$
there is an integer $n \geq m$ such that $A$ can be expressed
as $M_n(B)$ for a C*-algebra $B$ with $\tsr(B) \leq K$. 
Typical such an example is $B \otimes UHF$ for any unital \ca~$B$.
As a corollary, 
Let $A$ be a tsr boundedly divisible, unital C*-algebra with $\tsr(A) = 1$,
$G$ a finite group, and $\alpha$ an action of $G$ on $A$. Then
$\tsr(A \rtimes_\alpha G) \leq 2$.
This estimate is best possible. 
Indeed in \cite[Example 8.2.1]{bl3}
B. Blackadar constructed an symmetry action $\alpha$ on $CAR$
such that
$
(C[0,1] \otimes CAR) \rtimes_{id \otimes \alpha} Z_2 \cong C[0,1] \otimes B,
$
where $B$ is the Bunce-Deddens algebra of type $2^\infty$.
Then since $K_1(B)$ is non-trivial, we know that
$
\tsr(C[0,1] \otimes B) = 2.
$

In this note we try to solve generalized Blackadar's question and 
get the final estimate in some sense: 
Let $1 \in A \subset B$ be an inclusion of unital C*-algebras 
of index-finite type and depth $2$.
Suppose that $A$ is infinite dimensional simple with $\tsr(A) = 1$ 
and SP-property.
Then $\tsr(B) \leq 2$. 
In the case of crossed product algebras we conclude that 
$\tsr(A \rtimes_\alpha G) \leq 2$ for a simple unital 
C*-algebra $A$ with $\tsr(A) = 1$ and SP-property, and  an action $\af$ 
from   a finite group $G$ on $\Aut(A)$. 
We cannot still conclude that $\tsr(A \rtimes_\alpha G) = 1$, but it 
seems to guarantee that the question  would be solved affirmatively.

\section{Preliminaries}

\begin{dfn}\label{dfn2.1}
Let $A$ be a unital C*-algebra and $Lg_n(A)$ be
the set of elements $(b_i)$ of $A^n$ such that
$$
Ab_1 + Ab_2 + \cdots + Ab_n = A.
$$
 Then topological stable rank of $A$,
$\tsr(A)$, is defined to
be the least integer $n$ such that the set $Lg_n(A)$ is dense in $A^n$.
Topological stable rank of a non-unital C*-algebra is defined by
topological  stable rank of its unitaization algebra $\tilde{A}$
\end{dfn}

Note that $\tsr(A) = 1$ is equivalent to having the dense set of
invertible elements in $\tilde{A}$. 

The following is a  well-known characterization of topological stable rank one.
See \cite{Rf1} and \cite[Remark~2.4]{OT}.

\begin{prp}\label{P:Stablerankone}
Let $A$ be a unital \ca.
\begin{enumerate}
\item[$(1)$]
Let $p$ be a non-zero projection in $A$.
Then $\tsr(A) = 1$ if and only if $\tsr(pAp) = \tsr((1- p)A(1- p)) = 1$.
\item[$(2)$]
Let ${\mathbb K}$ be a C*-algebra of compact operators on an
 infinite dimensional separable Hilbert space. Then
$$
\tsr(A) = 1 \ \hbox{if and only if} \ \tsr(A \otimes {\mathbb K}) = 1.
$$
\end{enumerate}
\end{prp}

\begin{rmk}\label{R:Stablerankone}
Suppose that a \ca~$B$ is stably isomorphic to a \ca~$A$, that is,
$B \otimes {\mathbb K} \cong A \otimes {\mathbb K}$. 
Then from Proposition~$\ref{P:Stablerankone}(2)$ if $\tsr(A) = 1$, 
then $\tsr(B) = 1$.
\end{rmk}

Let $A$ be a \ca. $A$ is said to have {\it SP-property} if 
any non-zero hereditary C*-subalgebra of $A$ has non-zero 
projection. 
It is well known that if $A$ has real rank zero, that is, 
any self-adjoint element can be approximated by self-adjoint elements 
with finite spectra, then $A$ has SP-property. (See \cite{bp0}.)

\vskip 3mm

Next we summarize the C*-index theory of Watatani (\cite{wata}).

Let $1 \in A \subset B$ be an inclusion  of C*-algebras, and let $E\colon B
\rightarrow A$ be a faithful conditional expectation from $B$ to $A$.

A finite family
$\{(u_1, v_1), \dots, (u_n, v_n)\}$ in $B \times B$ is called {\em a

quasi-basis} for $E$ if $$
\sum_{i=1}^n u_iE(v_ib) = \sum_{i=1}^nE(bu_i)v_i = b \enskip \hbox{for}
\enskip b \in B.
$$
We say that a conditional expectation $E$ is  of {\em index-finite type}
if there
exists a quasi-basis for $E$. In this case the index of $E$ is defined by
$$
{\rm Index}(E) = \sum_{i=1}^nu_iv_i.
$$
(We say also  that the inclusion $1 \in A \subset B$ is
of {\em index-finite type}.)

Note that ${\rm Index}(E)$ does not depend on the choice of a quasi-basis
(\cite[Example~ 3.14]{iz3})
and every conditional expectation $E$ of index-finite type on a C*-algebra
has a quasi-basis of the form $\{(u_1, u_1^*), \dots, (u_n,u_n^*)\}$
(\cite[Lemma~ 2.1.6]{wata}).
Moreover ${\rm Index}(E)$ is always contained in the center of $B$, so
that it is a scalar whenever $B$ has the trivial center, in particular when
$B$ is simple (\cite[Proposition~ 2.3.4]{wata}).

Let $E\colon B\to A$ be a faithful conditional expectation.
Then $B_{A}(=B)$ is
 a pre-Hilbert module over $A$ with an $A$-valued inner
product $$\langle x,y\rangle =E(x^{*}y), \ \ x, y \in B_{A}.$$
Let $\mathcal E$ be the completion of $B_{A}$ with respect to the norm on
$B_{A}$ defined by
$$\| x\|_{B_{A}}=\|E(x^{*}x)\|_{A}^{1/2}, \ \ x \in B_{A}.$$
Then $\mathcal E$ is a Hilbert $C^{*}$-module over $A$.
Since $E$ is faithful, the canonical map $B\to \mathcal E$ is injective.
Let $L_{A}(\mathcal E)$ be the set of all (right) $A$-module homomorphisms
$T\colon \mathcal E\to \mathcal E$ with an adjoint $A$-module homomorphism
$T^{*}\colon \mathcal E\to \mathcal E$ such that $$\langle T\xi,\zeta
\rangle =
\langle \xi,T^{*}\zeta \rangle \ \ \ \xi, \zeta \in \mathcal E.$$
Then $L_{A}(\mathcal E)$ is a $C^{*}$-algebra with the operator norm
$\|T\|=\sup\{\|T\xi \|:\|\xi \|=1\}.$ There is an injective
$*$-homomorphism $\lambda \colon B\to L_{A}(\mathcal E)$ defined by
$$\lambda(b)x=bx$$
for $x\in B_{A}$ and  $b\in B$, so that $B$ can
be viewed as a
$C^{*}$-subalgebra of $L_{A}(\mathcal E)$.
Note that the map $e_{A}\colon B_{A}\to B_{A}$
defined by $$e_{A}x=E(x),\ \ x\in
B_{A}$$ is
bounded and thus it can be extended to a bounded linear operator, denoted
by $e_{A}$ again, on $\mathcal E$.
Then $e_{A}\in L_{A}({\mathcal E})$ and $e_{A}=e_{A}^{2}=e_{A}^{*}$; that
is, $e_{A}$ is a projection in $L_{A}(\mathcal E)$.
A projection $e_A$ is called the {\em Jones projection} of $E$.

The {\sl (reduced) $C^{*}$-basic construction} is a $C^{*}$-subalgebra of
$L_{A}(\mathcal E)$ defined to be
$$
C^{*}(B, e_{A}) = \overline{ span \{\lambda (x)e_{A} \lambda (y) \in
L_{A}({\mathcal E}): x, \ y \in B \ \} }^{\|\cdot \|} $$
(\cite[Definition~ 2.1.2]{wata}).

Since $C^{*}(B, e_{A})$ is isomorphic to $qM_n(A)q$ for 
some $n \in \N$ and projection $q \in M_n(A)$ by \cite[Lemma~3.3.4]{wata}, 
we conclude that
if $\tsr(A) = 1$, then $\tsr(C^{*}(B, e_{A})) = 1$ from 
Proposition~$\ref{P:Stablerankone}$ and Remark~$\ref{R:Stablerankone}$.

The inclusion $1 \in A \subset B$ of unital C*-algebras of
index-finite type is said to
have {\it finite depth k}
if the derived tower obtained by iterating the basic construction
$$
A' \cap A \subset A' \cap B \subset
A' \cap B_2 \subset A' \cap B_3 \subset
\cdots
$$
satisfies
$(A' \cap B_k)e_k(A' \cap B_k) = A' \cap B_{k+1}$,
where $\{e_k\}_{k \geq 1}$ are projections
derived obtained by iterating the basic construction such that
$B_{k+1} = C^*(B_{k}, e_k)$ \ ($k \geq 1$) \quad ($B_1 = B, e_1 = e_A$).
Let $E_k : B_{k+1} \rightarrow B_k$ be a faithful conditional expectation
correspondent to $e_
k$ for $k \geq 1$.

When $G$ is a finite group and $\alpha$ an
action of $G$ on $A$,
it is well known that an inclusion
$1 \in A \subset  A \rtimes_\alpha G$
is of depth 2. (See \cite[Lemma~3.1]{OT}.)

\section{Main result}

The following result is contained in \cite[Theorem~5.1]{OT}. 
We give a sketch of the proof for self-contained.

\begin{prp}\label{prpOT2}$($cf\cite[Theorem 5.1]{OT}$)$
Let $1 \in A \subset B$ be an inclusion of unital C*-algebra
of index-finite type and depth $2$. Suppose that $\tsr(A) = 1.$
Then we have
$$
\sup_{p\in P(A)}\tsr(pBp) < \infty,
$$
where $P(A)$ denotes the set of all prjections in $A$.
\end{prp}

\begin{proof}
Let
$$
1 \in A \subset B \subset B_2 \subset B_3 \subset \cdots
$$
be the derived tower of iterating the basic construction and
$\{e_k\}_{k\geq1}$ be canonical projections such that
$B_{k+1} = C^*(B_k,e_k)$, where $e_1 = e_A$.
Since $1 \in A \subset B$ is of depth 2,
we have
$$
(A' \cap B_2)e_2(A' \cap B_2) = A' \cap B_3.
$$
(See \cite[Definition 4.6.4]{GHJ}.)
Then
there are
finitely elements $\{u_i\}_{i=1}^n$ in  $A' \cap B_2$
such that
$$
\sum_{i=1}^nu_ie_2u_i^* = 1.
$$

Since for any $b \in B_2$

$$
\begin{array}{ll}
be_2 = 1\cdot be_2 &= \sum_{i=1}^nu_ie_2u_i^*be_B\\
&= \sum_{i=1}^nu_iE_2(u_i^*b)e_2,
\end{array}
$$
from \cite[Lemma~2.1.1 (2)]{wata}
we have  $b = \sum_{i=1}^nu_iE_2(u_i^*b)$. 
Similarly, we can show that 
$b = \sum_{i=1}^nE_2(bu_i)u_i^*$ for any $b \in B_2$.
It follows that
$\{(u_i, u_i^*)\}_{i=1}^n$ is a quasi-basis for $E_2$.

Since $pu_i = u_ip$ for $1 \leq i \leq n$ and any projection $p \in A$,
from the simple calculation we know that $\{(pu_i, pu_i^*)\}_{i=1}^n$
is a quasi-basis for $F_p = E_2|pB_2p$ from $pB_2p$ onto $pBp$.
Hence from \cite[Proposition~5.3]{OT} we have
$$
\begin{array}{ll}
\tsr(pBp) &\leq n^2\times\tsr(pB_2p) - n + 1\\
&= n^2 - n + 1.
\end{array}
$$
The last equality comes from that $\tsr(B_2) = 1$ and 
Proposition~$\ref{P:Stablerankone}(1)$.
Since $n$ is independent of the choice of projections in $A$,
we have 
$$
\sup_{p\in P(A)}\tsr(pBp) \leq n^2 - n + 1 < \infty.
$$
\end{proof}

The following main theorem is an extended version of 
\cite[Theorem~5.1]{OT}.

\begin{thm}\label{T:Main}
Let $1 \in A \subset B$ be an inclusion of unital C*-algebras 
of index-finite type and depth $2$.
Suppose that $A$ is infinite dimensional simple with $\tsr(A) = 1$ 
and SP-property.
Then $\tsr(B) \leq 2$.
\end{thm}

\begin{proof}
Since $1 \in A \subset B$ be an inclusion of simple C*-algebras 
of index-finite type and depth $2$, from Proposition~$\ref{prpOT2}$
we have 
$$
\sup_{p\in P(A)}\tsr(pBp) < \infty.
$$
Set $K = \sup_{p\in P(A)}\tsr(pBp)$. 

Since $A$ is simple with SP-property, there is a sequence of 
mutually orthogonal equivalent projections $\{p_i\}_{i=1}^N$
in $A$ such that $N > K$. (For example see \cite[Lemma 3.5.7]{Hl}.) 

Set $p = \sum_{i=1}^Np_i$. Then $pBp$ has a matrix unite such that
$$
pBp \cong M_N(p_1Bp_1).
$$
Then using \cite[Theorem 6.1]{Rf1} 
\begin{align*}
\tsr(pBp) &= \tsr(M_N(p_1Bp_1))\\
&= \{\frac{\tsr(p_1Bp_1) - 1}{N}\} + 1\\
&\leq \{\frac{K}{N}\} + 1 = 2,\\
\end{align*}
where $\{a \}$ denotes least integer greater than $a$. 
Since $A$ is simple, $p$ is a full projection in $A$, and moreover, in $B$. 
Hence from \cite[Theorem~4.5]{bl4} we have 
\begin{align*}
\tsr(B) \leq \tsr(pBp) \leq 2.
\end{align*}
\end{proof}

\begin{cor}
Let $A$ be a simple \ca~ 
with $\tsr(A) = 1$ and SP-property 
and $\af$ an action of a finite group $G$ on $\Aut(A)$.
Then 
$$
\tsr(A \rtimes_\af G) \leq 2.
$$
\end{cor}

\begin{proof}
If $A$ is finite dimensional, then $A \rtimes_\af G$ has finite dimensional,
and $\tsr(A \rtimes_\af G) = 1$. 
Hence we may assume that $A$ is infinite dimensional.

Since $A \subset A \rtimes_\af G$ is an inclusion of index-finite and 
depth 2 from \cite[Lemma~3.1]{OT}, it follows from Theorem~$\ref{T:Main}$.
\end{proof}

\begin{cor}\label{C:Crossed product}
Let $A$ be a simple \ca~ of tracial topological rank zero and 
$\af$ an action of a finite group $G$ on $\Aut(A)$.
Then 
$$
\tsr(A \rtimes_\af G) \leq 2.
$$
\end{cor}

\begin{proof}
Since $A$ has tracial topological rank zero, 
$A$ has $\tsr(A) = 1$ and SP-property. 
(For example see \cite[Lemma 3.6.6 and Theorem~3.6.10]{Hl}.)
Hence the conclusion comes from Corollary~$\ref{C:Crossed product}$.
\end{proof}

\begin{rmk}\label{blackadr's example}
If a given \ca~$A$ has only the condition 
of $\tsr(A) = 1$, the estimate in 
Corollary~$\ref{C:Crossed product}$
is best possible.
Indeed in \cite[Example 8.2.1]{bl3}
Blackadar constructed an symmetry action $\alpha$ on $CAR$
such that
$$
(C[0,1] \otimes CAR) \rtimes_{id \otimes \alpha} Z_2 \cong C[0,1] \otimes B,
$$
where $B$ is the Bunce-Deddens algebra of type $2^\infty$.
Then since $K_1(B)$ is non-trivial, we know that
$$
\tsr(C[0,1] \otimes B) = 2.
$$
$($See also \cite[Proposition~ 5.2]{nop}.$)$
\end{rmk}

\begin{rmk}
From Corollary~$\ref{C:Crossed product}$
if $A$ is infinite dimensional simple AF \ca~, we conclude that 
\begin{align*}
\tsr(A \rtimes_\af G) \leq 2
\end{align*}
for any an action $\af$ from any finite group $G$ on $Aut(A)$.
This gives an affirmative data  to  Blackadar's question in \cite{bl3},
that is, $\tsr(A \rtimes_\af G) = 1$ under the above condition. 
We note that under the extra condition on an action we can conclude 
it. Indeed, Phillips proved in \cite{Ph} if $\af$ has strictly Rokhlin 
property, then $A \rtimes_\af G$ is again AF \ca. More recently, 
the first author and Phillips proved in \cite{OP} that
if $\af$ is an action with tracial Roklin property from a finite group 
$G$ on a simple C*-algebra $A$ of tracial topological rank zero, then 
$A \rtimes_\af G$ has again tracial topological rank zero.
All data implies that Blackadar's question should be  correct. 
\end{rmk}

\begin{rmk}\label{R:Cancellation}
A \ca\ $A$ is said to have {\it cancellation of
projections} if   $p \sim q$ holds whenever $p, q$, and $r$ are
projections in $A$ with $p \perp r, \ q \perp r, \ p+r \sim q+r$.
If the matrix algebra $M_n(A)$ over $A$ has cancellation of
projections for each $n \in \N$, we simply say that $A$ has {\it
cancellation}.
It is well known that if $A$ has $\tsr(A) = 1$, then $A$ has cancellation.
But it was a long standing
problem till  quite recently whether the cancellation of $A$
implies $\tsr(A)=1$ for a stably finite simple \ca~ $A$, and it
has been settled  finally  in \cite{To} where a stably finite
simple $C^*$-algebra $B$ with cancellation and $\tsr(B)>1$ is
constructed by applying Villadson's techniques (\cite{v})  ($B$ is
also unital simple, separable, and nuclear).
Very recently, authors, Jeong, and Phillips have proved \cite{JOTP} that 
under the same assumption in Theorem~$\ref{T:Main}$ 
$B$ has cancellation. Therefore, we predict that 
$\tsr(B) = 1$.
\end{rmk}

\end{document}